\documentclass[11pt]{article}

\usepackage{amssymb,amsmath}
\usepackage{graphicx}

\newtheorem{theorem}{Theorem}[section]
\newtheorem{Proposition}{Theorem}[section]

\newtheorem{Lemma}[theorem]{Lemma}
\newtheorem{defi}[theorem]{Definition}

\numberwithin{equation}{section}

\newcommand{\non}{\nonumber}
\newcommand{\cR}{{\mathbb R}}

\newcommand{\cN}{{\mathbb N}}

\newcommand{\eq}[1]{\mbox{\rm {(\ref{#1})}}}

\usepackage{epsfig}

\title{\Large\bf Existence and regularity of weak solutions to a model\\
 for coarsening in molecular beam epitaxy}

\author{\small\sc Jun Zhang\thanks{E-mail:
     mathzj@zjut.edu.cn}\\
\small Department of Mathematics,\\
\small Zhejiang University of Technology\\
\small  Hangzhou 310032, Zhejiang\\
\small  P. R. China,\\
\small and\\
\small \sc Peicheng Zhu$^{1,2}$\thanks{E-mail: zhu@bcamath.org}\\
\small $^1$ Basque Center for Applied Mathematics (BCAM)\\
\small Building 500, Bizkaia Technology Park \\
\small E-48160 Derio, \ \ Spain\\
\small $^2$ IKERBASQUE, Basque Foundation for Science\\
\small E-48011 Bilbao,\ \ Spain}

\date{\ }

\begin{document}

\maketitle

\begin{abstract}

Taking into account the occurrence of a zero of the surface diffusion current and the requirement of
 the Ehrlich-Schwoebel effect, Siegert et al \cite{Siegert94} formulate a model of Langevin type that
 describes the growth of pyramidlike  structures on a surface under conditions of molecular beam epitaxy, and that
 the slope of these pyramids is selected by the crystalline symmetries of the growing film.
 In this article, the existence and uniqueness of weak solution to an initial boundary value problem
 for this model is proved, in the case that the noise is neglected.
 The regularity of the weak solution to models, with/without slope selection, is also investigated.



\medskip
\noindent{\bf  Running title:} Solutions to a model for coarsening dynamics.

\end{abstract}

\section{Introduction}

Many processes occur mainly at surfaces of materials, such processes include crystal growth, catalytic
reactions, production of nano-structures. Thus
surfaces are of great technological and fundamental interest.
There has been increasing interest in the understanding of the kinetics of surface growth
 processes, see e.g. \cite{Krug97, Barabasi95}.
 The growth of a crystalline film from a molecular or atomic beam is commonly referred to as molecular beam epitaxy (MBE),
  which is among the most fined methods for the
 growth of thin solid films and is of great importance for applied studies \cite{Johnson94}.
 MBE takes place in high vacuum or ultra high vacuum, for instance, some $10^{-8}$ Pa.
 The most important aspect of MBE is the slow deposition rate that is typically less than $1000$nm per hour,
 so MBE allows the films to  grow epitaxially. The slow deposition rates require proportionally
 better vacuum to achieve the same impurity levels as other deposition techniques. In turn, it is possible by using
 this technique to grow high-quality crystalline materials and form structures with very high precision in the
 vertical direction.

 The mathematical theory has been developed since Langevin who proposed the equation named after him.
 The  model by Kardar, Parisi abd Zhang (KPZ) \cite{Kardar86} describes well
 growth process such as the Eden process \cite{Eden58}, ballistic deposition \cite{Vold59},
   and growth of various restricted solid-on-solid models \cite{Plischke87}.
  Thus the KPZ model has been widely accepted as a model for the growth of crystals, and also has been extended to various cases.
 Let us mention especially the following work: In \cite{Johnson94} the authors consider the system that has potential barriers near
 step edges that suppress the diffusion of adatoms to a lower terrace. This effect is now commonly called
 Ehrlich-Schwoebel effect. They take into account the step-flow regime and instability, and
   proposed a continuum equation to  model the growth in MBE, which is valid only at the early time as
   long as the slopes are much smaller than $1$. However this is too restrictive to describe the unstable three-dimensional
    growth of real materials in the later time regime.  In order to take into account the occurrence of a zero in the surface
    diffusion current and the requirement of Ehrlich-Schwoebel effect,
    Seigert et al \cite{Siegert94} thus introduce a current for a structure with a cubic symmetry so that
 the current is changed to the one that has a zero differing from $0$ so that the model can be also applied
 to the regime when the slope is much greater than $1$. In this article, we shall study this model, also the model without slope
 selection will also be investigated. Due to the forming of steps, pyramidlike surfaces, etc. during the crystal
 growth, it is more natural to assume that the
 initial data is in $L^2(\Omega)$ or $H^1(\Omega)$ than in $H^2(\Omega)$.

\medskip
To formulate the model, we need some notations. Let $x=(x_1,x_2)\in\Omega\subset\cR^2$ be
 material point, $\Omega$ is an open bounded set with smooth boundary $\partial \Omega$. Let
 $t$ be the time variable. $Q_t = (0,t)\times \Omega$. $h=h(t,x)$ is the height that is
 measured in a co-moving frame of reference and describes the
local position of the moving surface.  $\nabla_x  h$ is the gradient of $h$,
 and $\Delta_x  = \partial^2/\partial x_1 ^2 + \partial^2/\partial x_2 ^2$ is the Laplacian.
For simplicity of notations, we shall use the following notations
$$
 \xi=(p, q) = \nabla_x h.
$$

 Then the equation turns out to be
\begin{eqnarray}
 \frac{\partial h}{\partial t}  + \nu\Delta_x ^2 h + {\rm div}_x \left( J(\nabla_x h) \right)  = 0 .
  \label{eq1}
\end{eqnarray}
which is satisfied in $Q_T$, where $T$ is a given positive number. And the boundary and initial conditions are
\begin{eqnarray}
 \frac{\partial}{\partial n}    h    = 0,  \ \left( \nu\nabla_x \Delta_x h + J(\nabla_x h)\right)\cdot n =0, & &  {\rm on}\ [0,T]\times\partial\Omega,
 \label{eq3a} \\
 h(0,x) = h_0(x), & &  x\in \bar\Omega.
 \label{eq3}
\end{eqnarray}
Here $n=(n_1,n_2)$ is the unit outward normal vector to the boundary $\partial \Omega$.
We have introduced the surface diffusion current $J = J(\nabla_x h) = J(p,q) $
\begin{eqnarray}
 J &=&  (j_{1},\ j_{2}),
 \label{flux0} \\
 j_{1} &=& \alpha \Big( (p+q)f\left((p+q)^2\right)  + (p-q)f\left((p -q)^2 \right) \Big),\\
 j_{2} &=& \alpha \Big( (p+q)f\left((p+q)^2\right)  - (p-q)f\left((p -q)^2 \right) \Big),
 \label{flux}
\end{eqnarray}
where $f$ is defined by
\begin{eqnarray}
 f(y) = \frac{1-y}{(1-y)^2 + \beta y},
 \label{nonlinearity}
\end{eqnarray}
and $ \alpha $ is a constant of surface diffusion, $\beta= (\ell_d)^2$ where $\ell_d$ is the diffusion
length.

This completes the formulation of an initial-boundary value problem.
 It is worth a remark on the nonlinearity $f$ since there are several varieties of $f$
 which lead to different models relating to ours.

\medskip
\noindent{\bf Remark 1.} {\it  Define
\begin{eqnarray}
 f(y) = \frac{1}{1 + \beta y} .
 \label{nonlinearityJohnson}
\end{eqnarray}
 Then the corresponding model is proposed by Johnson et al. {\rm \cite{Johnson94}}.
However this model does have a slope selection mechanism,  and  is correct only for
 early time as long as the slopes are much smaller than $1$. It is too restrictive to
 describe the unstable three-dimensional growth of real materials in the later time regime. Therefore, \eq{nonlinearity} is introduced
 by Siegert et al in {\rm \cite{Siegert94}} to interpolate the two regimes.
 The form of the surface current $J$ is the minimal model in the sense that the nonlinearity must be chosen
such that it describes the instability and leads to slope selection. The flux still
has the correct physical behavior: $|J|\sim \sqrt{ p^2+q^2 }$ for $p^2+q^2 \ll  \frac1{\ell_d}$ and
$|J|\sim 1/\sqrt{ p^2+q^2 }$ for $ \frac1{\ell_d}\ll p^2+q^2 \ll  1 $. This type of fluxes gives rise a completely different
behavior than the one defined by \eq{nonlinearityJohnson}, despite many similarities, as shown in {\rm \cite{Siegert94}}.
 The exact form of $f$ does not
 play a role since the slope selection mechanism and the growth exponents do not depend on such details.

In {\rm \cite{Stroscio95,Rost97}}, the current of the form
\begin{eqnarray}
 J = \xi (1-|\xi|^2)
 \label{nonlinearity1}
\end{eqnarray}
is used, however, it has stable zeros for all slopes with $|\xi|=1$ regardless of the direction of $\xi$. Thus such an azimuthal
 symmetry is unrealistic for crystalline films.
Therefore here $f_i$ ($i=1,2$) are functions chosen such that $f_1(p^2,q^2) = f_2(q^2,p^2)$.
The simplest form that describes growth on such substrates is a current with components
\begin{eqnarray}
 j_1 &=& p (1- p^2 - b q^2),
 \label{nonlinearity2a} \\
 j_2 &=& q (1- q^2 - b p^2),
 \label{nonlinearity2}
\end{eqnarray}
which leads to a buildup of pyramids with selected slopes $(p_0,q_0) = (\pm1,\pm1)/\sqrt{1+b} $ for $-1< b <1$. This diffusion
current is suitable for substrates with a quadratic symmetry.

Finally we point out that after a coordinate transformation $X=Ax$,
where
$$
 A=\left(\begin{array}{cc}
 1 & 1 \\
 1 & -1
\end{array}
\right),
$$
we find $J$ in \eq{nonlinearity} can be reduced, without loss of generality, to a simpler form
\begin{eqnarray}
 j_{1} &=& \alpha \, p f(p^2), \\
 j_{2} &=& \alpha \,  q f(q^2).
 \label{flux1}
\end{eqnarray}

}

It is interesting to compare the important difference between the Cahn-Hilliard equation modeling
 phase-ordering and the model considered here.
 %
{Mathematically,  equation \eq{eq1}, with nonlinearity \eq{nonlinearityJohnson}, \eq{nonlinearity1}
   or \eq{flux1}, differs from the Cahn-Hilliard equation due to the flux term: $J$ in this paper depends on
the gradient of the unknown, while the  Cahn-Hilliard equation on the unknown only. Numerical experiments also
show the important differences between this model and the Cahn-Hilliard one.  Many papers, e.g. \cite{Siegert94,Siegert98,Siegert97} have
 been carried out the study of the differences between the problem studied here and phase-ordering
dynamics described by the Cahn-Hilliard equation \cite{Cahn58}. These differences become apparent when the domain configurations are
plotted as in Figure 1. A domain in this context is an area of constant slope corresponding to one of the four values.
 The analogous case in phase-ordering dynamics is described by a four-state clock model, see, e.g. \cite{Siegert97, Foret05}.
  However, in that case we shall find that domain walls do not have any particular orientation,
whereas here domain walls are intersections of planes of
constant slopes and therefore form {\it straight lines}. Furthermore, there are two types of domain walls: Domain walls
at which only one component of the slope changes are
aligned along the $x$ and $y$ axes. These are the edges of the
pyramids;   Domain walls at which both components of the slope change run at $45^\circ$ with respect to the
principal axes. These latter domain walls form the so-called roof tops.

%
\vskip3.68cm
 \includegraphics[scale=0.5]{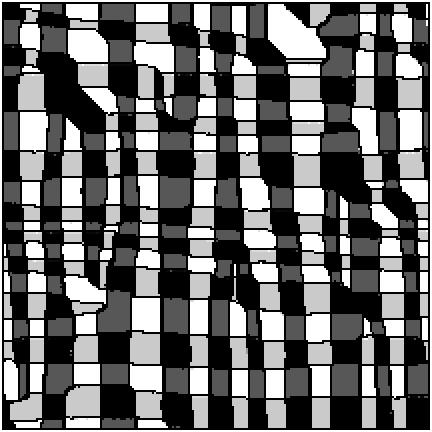}
%
%
 \hskip6.05cm
 \includegraphics[scale=0.96]{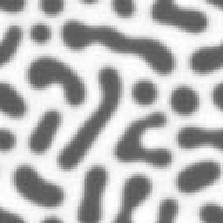}

\noindent Figure 1. Left: the configuration for the crystal growth model, from \cite{Siegert98};
Right: the domain walls of phase-ordering governed by the Cahn-Hilliard equation, from e.g. \cite{Foret05}.

}

\bigskip
Now we derive the model briefly. Define the free energy by
$$
 E[h] = \int_\Omega \left( \frac\nu2 (\Delta_x  h)^2 - \frac{\alpha}2 \left(F\left( (h_{x_1}+h_{x_2})^2\right)
 + F\left((h_{x_1} - h_{x_2})^2\right)\right)   \right)dx.
$$
Here $F = \int^x f(y)dy$ is the primitive of $f$ such that $F(x)$ grows at $-\log(|x|)$ as $|x|$ goes
to $\infty$. Suppose now that $h=h(t,x)$ is a solution to \eq{eq1} -- \eq{eq3}. Formal computations yield
\begin{eqnarray}
 & &\frac{d}{dt} E[h] =  \int_\Omega \Big( - \alpha f\left( ( h_{x_1} + h_{x_2})^2\right)(h_{x_1}+h_{x_2})(h_{x_1} +h_{x_2})_t  \Big)dx\non\\
 & & + \int_\Omega \Big(
 - \alpha  f\left( ( h_{x_1} - h_{x_2})^2 \right) (h_{x_1}-h_{x_2})(h_{x_1} - h_{x_2})_t  +  \nu   \Delta_x  h \Delta_x  h_t   \Big)dx.\non
\end{eqnarray}
Combining the terms containing $(h_{x_1})_t$ (or $(h_{x_2})_t$) together we can rewrite the right hand side of the above equality as
\begin{eqnarray}
 &&  -\int_\Omega   \alpha \Big(  f\big( ( h_{x_1} + h_{x_2})^2\big)(h_{x_1}+h_{x_2})
 +  f\left( ( h_{x_1} - h_{x_2})^2\right)(h_{x_1}-h_{x_2})  \Big) (h_{x_1}  )_t\, dx  \non\\
 & &   - \int_\Omega \alpha \Big( f\left( ( h_{x_1} + h_{x_2})^2\right)(h_{x_1}+h_{x_2})
 -  f\left( ( h_{x_1} - h_{x_2})^2\right)(h_{x_1}-h_{x_2})  \Big) (h_{x_2}  )_t \,   dx \non\\
 & & + \int_\Omega  \nu \Delta_x  h \Delta_x  h_t \,  dx\non\\
 &=&  \int_\Omega \left( - J\cdot \nabla_x h _t        +  \nu   \Delta_x  h \Delta_x  h_t   \right)dx.
 \label{secondlaw1}
\end{eqnarray}
Using integration by parts and equation \eq{eq1}, we infer from \eq{secondlaw1} that
 \begin{eqnarray}
  \frac{d}{dt} E[h]  &=&    \int_\Omega  \left( - J - \nu  \nabla_x \Delta_x  h     \right)\cdot \nabla_x h_t\,  dx\non\\
  &=&   \int_\Omega {\rm div}_x \left( J  + \nu  \nabla_x \Delta_x  h     \right) h_t\,  dx\non\\
 &=& - \int_\Omega \left( {\rm div}_x \left( J   + \nu  \nabla_x \Delta_x  h    \right) \right)^2 dx\non\\
 &\le& 0.
 \label{secondlaw}
\end{eqnarray}
This implies the second law of thermodynamics is valid. The equation can be written in a gradient form with the total surface
 current defined by $J_1 :=  J_{eq} + J(\nabla_x h)$ where $J_{eq}:=\nu  \nabla_x \Delta_x  h$ is called the equilibrium surface current,
 $\nu$ is proportional to the surface stiffness.

\medskip
\noindent{\bf Statement of the main result.} Before the statement of our main results, we define weak solutions to problem
consisting of (\ref{eq1}) --  (\ref{eq3}). We use the notations: $(f,g)_{Q_{T}}$  and $(f,g)_{\Omega}$ are, respectively, the inner
 product of $f$ and $g$ over $Q_{T}$ and $\Omega$.  $\langle f,g\rangle_{Q_{T}}$ denotes the dual product of
$ f,g$ with $f\in L^2(0,T;X'),\ g \in L^2(0,T;X )$ and $X$ is a Banach space and $X'$ is its dual.
 $H^m(\Omega)$ are the standard Hilbert spaces of order $m$.
 Define $H_N^2(\Omega):=\{f\in H^2(\Omega)\mid \frac{\partial f}{\partial n} = 0\ {\rm on}\
 \partial \Omega\}$ and its dual space is denoted by $H_N^{-2}(\Omega)$.
%
\begin{defi} \label{D1.1}
Let $h_0\in L^2(\Omega)$. A function $h=h(t,x)$ with
\begin{eqnarray}
h \in L^\infty(0,{T};L^2(\Omega)) \cap L^2(0,T;H^2(\Omega)),\ h_t\in L^2(0,T;H_N^{-2}(\Omega))
\label{property}
\end{eqnarray}
is a weak solution of  problem (\ref{eq1}) -- (\ref{eq3}), if
\begin{eqnarray}
 \langle h_t,\varphi\rangle _{Q_{T}} + \nu ( \Delta_x h ,\Delta_x \varphi )_{Q_{T}}
  + \left(J , \nabla_x\varphi \right)_{Q_{T}}   = 0
\label{definition}
\end{eqnarray}
holds for all $\varphi\in L^2 (0,T; H_N^2(\Omega))$, and
$\lim_{t\downarrow 0^+}( h(t), \psi )_\Omega = (h_0,\psi )_\Omega$ for all $\psi\in L^2(\Omega)$.

\end{defi}
%

Now we are in a position to state the main results of this article.
%
%
\begin{theorem}[Existence] \label{T1.1} Suppose that the boundary of\, $\Omega$ is smooth,
and $h_0\in L^2(\Omega)$.  Then there exists a unique weak
solution $h$ of problem \eq{eq1} --  \eq{eq3} in the sense of Definition~1.1,
and  the total mass $\int_\Omega h(t,x)dx$ is conserved, i.e.
 $\int_\Omega h(t,x)dx\equiv \int_\Omega h_0(x)dx$.

Moreover, if $h_0\in H^1(\Omega)$, the weak solution of problem \eq{eq1} --  \eq{eq3} has, which in addition to
(\ref{property}), the following regularities
\begin{equation}
 h \in L^\infty(0,T; H^{1}(\Omega)) \cap L^2(0,T; H^{3}(\Omega)),\ h_t\in L^2(0,T; H^{-1}(\Omega)).
 \label{proper1}
\end{equation}

\end{theorem}

\begin{theorem}[Regularity]\label{T1.2} Suppose that  $h_0\in H^{2m}(\Omega)$ with $m\in\cN$.  Then the weak
solution $h$ of problem \eq{eq1} --  \eq{eq3} satisfies
\begin{eqnarray}
 h \in L^\infty(0,T; H^{2m}(\Omega)) \cap  L^2(0,T; H^{2m+2}(\Omega)),\ h_t \in  L^2(0,T; H^{2m-2}(\Omega)),
 \label{proper2a}
\end{eqnarray}
and
\begin{eqnarray}
 \begin{array}{ll}
  D^{l}_t h\in L^\infty(0,T; L^{2}(\Omega)) \cap L^2(0,T; H^{2}(\Omega)), & if\ m=2l, \ l\in\cN;\\[0.2cm]
  D^{l-1}_t h\in L^\infty(0,T; H^{2}(\Omega)) ,\ D^{l}_t h\in  L^2(0,T; L^{2}(\Omega)) & if\ m=2l-1, \ l\in\cN.
  \end{array}
 \label{proper2}
\end{eqnarray}
Consequently, if $ h_0\in C^{\infty}(\bar\Omega)$, then the solution $h$ is smooth on $\bar{Q}_{T}$.
\end{theorem}

\medskip
Now let us recall some references related closely to our problem. In \cite{King03}  an initial boundary value
 problem of this epitaxial model with cubic nonlinearities is studied in which
the initial data is chosen in $H^2(\Omega)$, and the $H^2(\Omega)$-norm of the solution follows directly from
the Clausius-Duhem inequality, the second laws of thermodynamics, but this technique does not work for our case since
we assume the initial data is only in $H^1(\Omega)$. Li and Liu study the initial boundary value problem for the MBE model
 with or without slope selection in \cite{Li} and
   the boundary conditions   are chosen periodic. In both articles, they construct approximate solutions by using the Galerkin method,
   while we use a linearized problem, together with the convolution technique, to obtain a sequence of smooth approximate solutions,
   then establish   a priori estimates for this sequence.
 Kohn and Otto \cite{Kohn} investigate the coarsening Rate for the Cahn-Hilliard equation, and Kohn and Yan \cite{Kohn03} studies
 the coarsening rate for an epitaxial growth model with a cubic nonlinearity. Watson and Norris \cite{Watson06} study the
 coarsening dynamics of multiscale solutions to a dissipative singularly perturbed partial differential equation with a trigonally
  symmetric potential which  models the evolution of a thermodynamically unstable crystalline surface.

We give now a remark on the choice of the initial data in this article.\\
\noindent{\bf Remark 2.} {\it The assumption that initial data $h_0$ is in $H^1$ is more natural than the one that $h_0$
is in $H^2$. The reason is that $h$ is piecewise affine in the case that the surfaces are high-symmetric, such as pyramidlike ones.
One evidence can be also seen from a typical Scanning Tunneling Microscope (STM) picture, see Figure~2,
 which shows clearly that the surface are not smooth. Correspondingly, a good mathematical model should consider this feature.

If the initial data is in  $H^1 $, the existence of  weak solutions to problem \eq{eq1} -- \eq{eq3} in which the nonlinearity is cubic,
 like \eq{nonlinearity1}, also \eq{nonlinearity2a} -- \eq{nonlinearity2}, is still open, and may be interesting.
  To solve such a problem, I surmise we need to invent an inequality of the Brezis-Gallouet type {\rm \cite{BG80}}.
\vskip.26cm

\begin{figure}[h!]
\hskip3.005cm
 \includegraphics[bb=10 20 100 300,angle=-90,width=0.204\textwidth]
{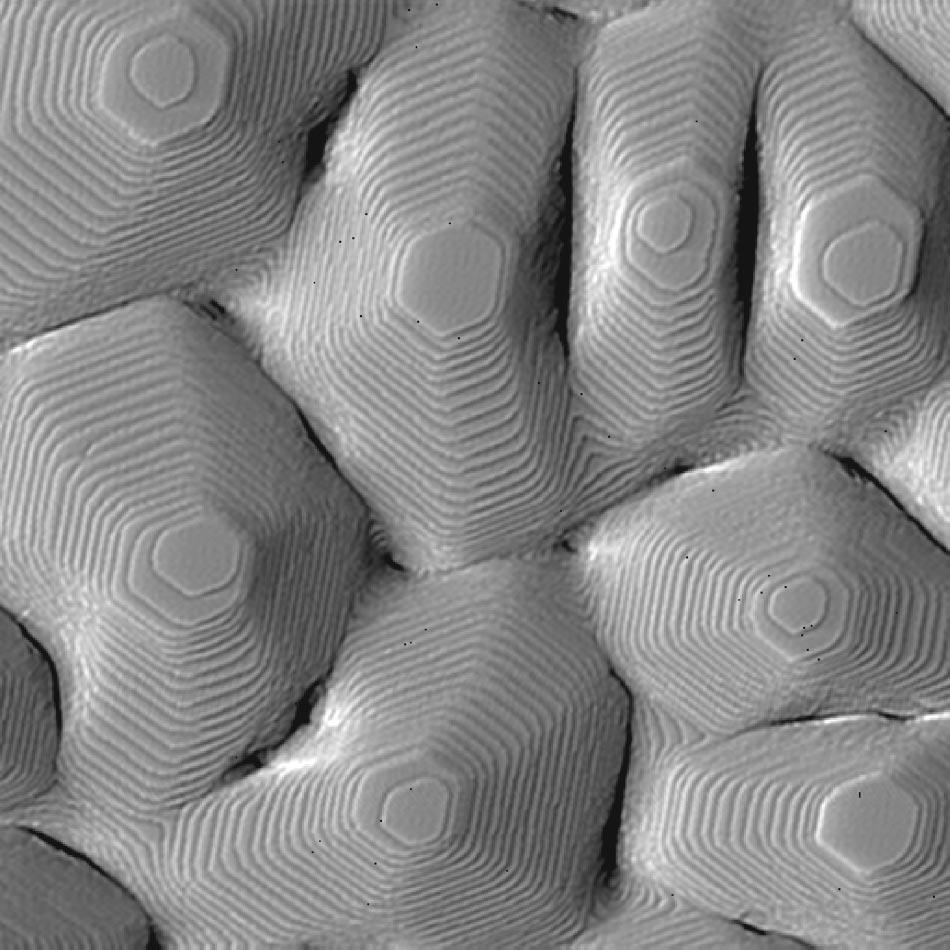}

\vskip6.668cm
\end{figure}
\noindent Figure 2. Typical STM image of mounds appearing on the surface.  From Krug, Politi and Michley {\rm \cite{Krug00}}.

\noindent
 }

The organization of the remaining parts of this article is as follows. The existence and uniqueness
 of weak solutions is studied in Section~2 by constructing smooth approximate solutions and using a priori estimates. Section~3 consists of
  two parts, each of which is concerned, respectively, with the regularities of weak solutions to models with and without slope selection.

%
\section{Existence of weak solutions}

\subsection{Existence for the approximate problem}

To prove the existence Theorem~\ref{T1.1}, we construct smooth approximate solutions in a similar way as done in \cite{Alber05},
 provided the initial data is smooth.
 Then we establish uniform a priori estimates of these solutions by which we conclude compactness. Before formulating
 an approximate problem, we introduce the modifier $\chi = \chi(t,x)$ such that
$
\chi\in C_0^\infty(Q_T)
$
satisfies
$
\int_{ \cR^3 } \chi(t,x)dtdx=1.
$
For $\varepsilon>0$, we set
$$
\chi_\varepsilon(t,x):=\frac1\varepsilon\chi\left(\frac{t}{\varepsilon}, \frac{x}{\varepsilon} \right),
$$
and for any function $f \in L^\infty(Q_{T})$ we define
\begin{eqnarray}
 \widetilde{f}(t,x) = (\chi_\varepsilon*f)(t,x)= \int_{ \cR^3 } \chi_\varepsilon (t-s,x-y) f(s,y) ds dy.
 \label{convolution}
\end{eqnarray}
We choose a smooth sequence $h_0^\varepsilon$ such that
 $$
 \|h_0^\varepsilon - h_0 \|_X\to 0
 $$
as $\varepsilon\to 0$. Here $X=L^2(\Omega)$ or $X=H^1(\Omega)$.

Then the smoothed problem turns out to be
\begin{eqnarray}
\frac{\partial h}{\partial t}  + \nu\Delta_x ^2 h + {\rm div}_x  ( J(\widetilde{\nabla_x \hat h})  )  = 0 ,
  \label{eq1Smoothed}
\end{eqnarray}
with the following boundary and initial conditions
\begin{eqnarray}
 \frac{\partial}{\partial n} h = 0,  \ \left( \nu\nabla_x \Delta_x h + J(\widetilde{\nabla_x \hat h}) \right) \cdot n =0, & &  {\rm on}\ [0,T]\times\partial\Omega,
 \label{eq3aSmoothed} \\
 h(0,x) = h_0^\varepsilon (x), & &  x\in \bar\Omega.
 \label{eq3Smoothed}
\end{eqnarray}

Note that equation \eq{eq1Smoothed} is a linear fourth order parabolic equation with a smooth known term.
By the existence theorem for higher order parabolic equations in the book by Ladysenskaya et al. \cite{Ladyzenskaya}
or the book by Eidelman \cite{Eidelman69},  we have
\begin{Proposition}[Existence of smooth approximate solutions] \label{T1.1a}
 Suppose that the assumptions of Theorem~\ref{T1.1} are satisfied, and $\varepsilon$ is a given
 positive constant. Let $\hat h\in L^2(0,T;H^2(\Omega))$.

Then for any $T>0$, there exists a unique smooth
solution $h$ of problem \eq{eq1Smoothed} --  \eq{eq3Smoothed}, which satisfies that
 the total mass $\int_\Omega h(t,x)dx$ is conserved.

\end{Proposition}

The solution $h$ constructed in Theorem~2.1 depends on the small parameter $\varepsilon$. In order to prove the existence
of weak solutions to the original problem \eq{eq1} --  \eq{eq3Smoothed}, we need to establish some a priori estimates which are
independent of $\hat h$ and $\varepsilon$ and thus guarantee
the passage to limit of $h^\varepsilon$ as $\varepsilon\to 0$.

\subsection{A priori estimates }

Assume that there exists a classical solution $h^\varepsilon$ to problem \eq{eq1Smoothed} -- \eq{eq3Smoothed}
 with smooth initial data $h_0^\varepsilon$ and $\hat h\in L^2(0,T;H^2(\Omega))$ satisfying
 $\|\hat h\|_{L^2(0,T;H^2(\Omega))}\le \bar C$. In what follows,  $C$ denotes a  constant which is independent of $\varepsilon$
 and $\hat h$. The $L^2(\Omega)$-norm of $f$ is denoted by $\|f\|$.
 We shall derive  a priori estimates for this solution. To begin with, we first state the following lemma on the nonlinearity $f$.
\begin{Lemma}\label{L2.0} There hold, for all $y\ge 0$, that
\begin{eqnarray}
 | y f(y) |  & \le & C,
  \label{est0a}\\
  | y^{m+1} f^{(m)}(y) |  & \le & C,\ m\in \cN.
  \label{est0b}
\end{eqnarray}
Here $f^{(m)} $ is the $m$-th order derivative of $f$.

\end{Lemma}

\noindent{\it Proof.}  One needs to investigate the behavior of the function
$$
 -yf(y) =  \frac{y^2 - y}{(1-y)^2 + \beta y}
$$
for all $y\ge 0$. It is easy to see that $-yf(y)\to 1$ as $y\to\infty$, hence $ |-yf(y)| \le C$ for $y\ge M$,
where $M$ is a suitably large constant. We shall prove that $-yf(y)$ is also bounded on the interval $[0,M]$. To this end,
we need to prove that the denominator $g(y) = (1-y)^2 +  \beta y$ is greater than a positive constant. Note
that $g(y)$ is nonnegative for $y\ge 0$. Thus $0$ is the only possible minimum of $g(y)$. However $g(y)=0$
implies that $1-y=0$ and $y=0$ which cannot be satisfied simultaneously. Therefore, the minimum of $g$ must be positive, and
since $M$ is finite, we infer from the continuity of $g(y)$    that
\begin{eqnarray}
 \min_{0\le y \le M} g(y)=C_1 >0 .
 \label{positivity}
\end{eqnarray}
We  then obtain, for all $y\in [0,M]$, that
$$
 |yf(y)|\le \frac{M^2 + M}{C}\le C.
$$

Next we consider the behavior of derivatives of $f$. Let $ m\le 1$ be an integer.
 Rewrite
$$
 f^{(m)}(y) = ((1-y) g(Y(y)) )^{(m)}
$$
where $g(Y) = Y^{-1}$ and $ Y = (1-y)^2 +  \beta y $.   Invoking the product rule
$$
f_0^{(m)} = (f_1\cdot f_2) ^{(m)}
  = \sum_{k=0}^m C_m^k  f_1 ^{(k)} f_2 ^{(m-k)} ,
$$
 where $C_m^k$ denotes the number of $k$-combinations of an $m$-element set, we have
\begin{eqnarray}
 f^{(m)}(y) &=& C_m^0(1-y) (g(Y(y)) )^{(m)} + C_m^1(1-y)' ( g(Y(y)) )^{(m-1)} \non\\
 &=&  (1-y)  (g(Y(y)) )^{(m)} - m (g(Y(y)) )^{(m-1)}.
 \label{productrule}
\end{eqnarray}
Making use of the Fa\`a di Bruno formula, i.e.
\begin{eqnarray}
 \frac{d^m}{dx^m} g(Y(y)) = \sum_{1 l_1+2l_2+\cdots  + m l_m = m}
 \frac{m!}{l_1!l_2!\cdots l_m!} g^{(l_1+l_2+\cdots +l_m)}(Y(y)) \prod
 _{j=1}^m\left(\frac{Y^{(j)} (y)}{j!}\right)^{l_j}
 \label{brunorule}
\end{eqnarray}
and recalling  $( g(Y) )^{(m)} = ( Y^{-1} )^{(m)} =(-1)^m m!\, Y^{-m-1} $,
$Y'(y)=-2(1-y) + \beta $, $Y''(y)=~2$ and $Y^{(j)} (y) = 0$ for any
$j>2$, we can reduce \eq{brunorule} to
\begin{eqnarray}
 \frac{d^m}{dx^m} g(Y(y)) = \sum_{1 l_1 + 2l_2 = m}(-1)^{l_1 + l_2}
 \frac{m!(l_1 + l_2)!}{l_1!\, l_2! }\, Y^{-(l_1+l_2+1)} (y) \prod_{j=1}^2 \left( \frac{Y^{(j)} (y)} {j!} \right)^{l_j}.
 \label{brunorule1}
\end{eqnarray}

For $y\in [0,M]$, $f^{(m)}(y)$ is smooth by \eq{positivity}. One thus needs only to
investigate the behavior for large $y\in [M,\infty)$, and it is enough to calculate
the highest exponent. From \eq{brunorule1} it follows
that the highest exponent is less than or equal to
$$
 {-2(l_1+l_2+1) + 1 l_1 + 0 l_2} =  {- l_1- 2l_2 -2} =  {- m -2 }.
$$
Therefore, invoking \eq{productrule}, we assert that there exists a constant $\gamma$ such that
\begin{eqnarray}
 f^{(m)}(y) \sim \gamma y^{- m - 1 }  \ {\rm as}\ y\to\infty.
 \label{productrule1}
\end{eqnarray}
Hence,  this implies \eq{est0b}. Thus the proof  of this lemma is complete.

\medskip
From now on we are going to derive a priori estimates. The first is
\begin{Lemma}\label{L2.1} There hold for all $t\in [0,T]$
\begin{eqnarray}
 \int_\Omega h(t,x) dx & = & \int_\Omega h_0^\varepsilon(x) dx,
  \label{est1a}\\
 \|h(t)\|^2 + \int_0^t  \|\Delta_x h(\tau)\|^2  d\tau &\le & C.
  \label{est1b}
\end{eqnarray}

\end{Lemma}

\noindent{\it Proof.}  Integrating \eq{eq1Smoothed} with respect to $x$ yields
\begin{eqnarray}
 \frac{d}{dt} \int_\Omega h (t,x) dx   = 0,
 \label{est1a1}
\end{eqnarray}
which implies \eq{est1a}.

Multiplying \eq{eq1Smoothed} by $h$ and integrating the resulting equation with respect to $x$, we obtain
\begin{eqnarray}
 \frac12\frac{d}{dt}\|h\|^2 + \nu\|\Delta_x h\|^2 + \int_\Omega J (\widetilde{\nabla_x \hat h}) \cdot \nabla_x h dx  = 0.
 \label{est1a2}
\end{eqnarray}
Applying Lemma~\ref{L2.0} we obtain
\begin{eqnarray}
 \int_\Omega J (\widetilde{\nabla_x \hat h}) \cdot \nabla_x h dx \le C \|\nabla_x h\|.
 \label{Wq}
\end{eqnarray}
Using the Poincar\'e inequality of the form $\|\nabla_x h\|\le C\|D_x^2 h\| +C(\int_\Omega h dx)$, applying the elliptic estimates
 $\|D_x^2 h\| \le C\|\Delta_x h\|$ and  from \eq{est1a}
 it then follows that
$$
 \|\nabla_x h\|\le  C\|\Delta_x h\| + C.
$$
Therefore, \eq{est1a2}  becomes
\begin{eqnarray}
 \frac12\frac{d}{dt}\|h\|^2 + \nu \|\Delta_x h\|^2 \le \frac\nu2 \|\Delta_x h\|^2  + C.
 \label{est1a3}
\end{eqnarray}
From this,  estimate \eq{est1b} follows. And the proof of this lemma is complete.

\begin{Lemma}\label{L2.2} There holds
\begin{eqnarray}
 \|h(t)\|^2_{H^1(\Omega)}+ \int_0^t \|\nabla_x \Delta_x h(\tau)\| ^2d\tau &\le & C.
  \label{est2a}
\end{eqnarray}

\end{Lemma}

\noindent{\it Proof.} Integrating, with respect to $x$  over $\Omega$, equation \eq{eq1Smoothed}  multiplied by $-\Delta_x h$ and  using
integration  by parts yield
\begin{eqnarray}
 \frac12\frac{d}{dt}\|\nabla_x h\|^2 + \nu\|\nabla_x \Delta_x h\|^2
 - \int_\Omega J(\widetilde{\nabla_x \hat h})\cdot \nabla_x  \Delta_x h\, dx  = 0.
 \label{est2a1}
\end{eqnarray}
Making use of   Lemma~\ref{L2.0} again, we can easily prove that
\begin{eqnarray}
   \left|\int_\Omega J(\widetilde{\nabla_x \hat h}) \cdot \nabla_x \Delta_x h dx \right|  \le   C +\frac\nu2 \|\nabla_x \Delta_x h\|^2   .
 \label{est2a2}
\end{eqnarray}

Thus integrating \eq{est2a1} with respect to $t$ one gets
\begin{eqnarray}
 \frac12 \|\nabla_x h\|^2 + \nu\int_0^t \|\nabla_x \Delta_x h\|^2d\tau \le C +  \frac\nu 2 \int_0^t \|\nabla_x \Delta_x h\|^2 d\tau.
 \label{est2a3}
\end{eqnarray}
From which estimate \eq{est2a} follows. The proof of this lemma is complete.

\subsection{Existence of solutions to the phase field model}

In this section we shall make use of the {\it a priori} estimates established in the previous
subsection to study the convergence of the  solutions $h^m $  of the
approximate problem for $m \to\infty$, thereby proving Theorem \ref{T1.1}.
In our investigation we need the following  well known results, see, for instance, Lions~\cite{Lions}, Evans~\cite{Evans90}:

We shall make use of the following lemma which is of Aubin-Lions type.
\begin{Lemma} \label{L2.6a}
Let $B_0,\ B,\ B_1$ be Banach spaces which satisfy
that $B_0,\ B_1$ are reflexive and that
$$
 B_0\subset\subset B\subset B_1.
$$
Here, by $\subset\subset$ we denote the compact imbedding.
Define
$$
 W=\left\{f \mid f\in L^\infty(0,T;B_0), \quad \frac{df}{dt}\in L^r(0,T;B_1) \right\}
$$
with $T$ being a given positive number and $1< r <\infty$.

Then the embedding of $W$ in $C([0,T]; B)$ is compact.

\end{Lemma}
%
To deal with the nonlinear terms, we also need
\begin{Lemma}\label{L5.2} Let $\Gamma$  be an open set in $\cR^m$.  Suppose
functions $f_n, f$ are in $L^q(\Gamma)$ for any given  $1<q<\infty$, which satisfy
$$
\|f_n\|_{L^q(\Gamma )}\le C, \ f_n\to f \ almost\ everywhere \
in\ \Gamma .
$$
Then $f_n$ converges to $f$ weakly in $L^q(\Gamma )$.
\end{Lemma}

\medskip
We now turn to prove the existence of weak solutions for the initial $h_0$ that is assumed in $L^2(\Omega)$. The existence
for the initial data in $H^1(\Omega)$ is easy by recalling Lemma~\ref{L2.2}.
  Using problem \eq{eq1Smoothed} -- \eq{eq3Smoothed}, we can construct smooth approximate
solutions as follows:  Let $h^1$ be a given function.  By solving problem \eq{eq1Smoothed} -- \eq{eq3Smoothed} with
$\hat h=h^1$, one gets $h^2$. Suppose we have obtained $h^m$  for some $m \in\cN$, then set
$\varepsilon = \frac1m$ and $\hat h = h^m$. Thus, we can define $h^{m+1}$ successively. Therefore we has a sequence of approximate solutions $h^m$.
 Since equation (\ref{eq1}) is nonlinear, we need some
 results about strong and pointwise convergence.

Define $f=h^j $ and $r=2$. Set
$$
 B_0=H^2(\Omega), \ B=H^1(\Omega), \ B_1=H^{-2}(\Omega),
$$
it follows from Lemma~\ref{L2.1} that
$$
 \|h\|_{L^2(0,T;B_0)} \le C,\ \|h_t\|_{L^2(0,T;B_1)}\le C.
$$
Applying Lemma~\ref{L2.6a}, we then conclude that $\{h^m\}$ is compact
in $C([0,T];B)$, namely  $C([0,T];H^1(\Omega))$, so $\{h^m_{x_i}\}$, for $i=1,2$, is compact in $C([0,T];L^2(\Omega))$.
Therefore there exists a subsequence, still denote it by $\{h^m_{x_i}\}$, such that
\begin{eqnarray}
 \|h^m - h\|_{C([0,T];H^1(\Omega))}\to 0,\ \ \|h^m_{x_i} - h_{x_i}\|_{C([0,T];L^2(\Omega))}\to 0
  \label{convergence0}
\end{eqnarray}
as $m\to\infty$. Moreover, we can select furthermore a subsequence such that
$ h^m_{x_i} $ converges to $h_{x_i}$ almost everywhere. Setting $\kappa=\frac1m$,  $S^\kappa = \nabla_x h^m$. By the properties of convolution
and \eq{convergence0}, we have
\begin{eqnarray}
 \|\chi_\kappa*S^\kappa-S \|_{{L^2}(Q_{T_e})}&\le &
 \|\chi_\kappa*(S^\kappa-S)\|_{{L^2}(Q_{T_e})}+
 \|(S-\chi_\kappa*S) \|_{{L^2}(Q_{T_e})}\non\\
 & \le & \|(S-\chi_\kappa*S) \|_{{L^2}(Q_{T_e})}
 + \|S^\kappa-S \|_{{L^2}(Q_{T_e})}\to 0,
 \label{convergConvolution}
\end{eqnarray}
for $\kappa\to 0$, whence we can select a subsequence, still denote by $\chi_\kappa*S^\kappa$ converges to $S$ almost everywhere.
Consequently,  we assert that
$$
 J(\widetilde{\nabla_x h^m})\ \text{converges to } J(\nabla_x h)
$$
almost everywhere as $m\to\infty$.  Remembering that $\|J(\widetilde{\nabla_x h^m})\|_{ L^2(Q_T)}\le C$, using Lemma~\ref{L5.2}
 we assert that
\begin{eqnarray}
  J(\widetilde{\nabla_x h^m})\rightharpoonup J(\nabla_x h)
  \label{convergence1}
\end{eqnarray}
in $L^2(Q_T)$ as $m\to\infty$.

For the linear terms, by weak compactness, one can easily get
\begin{eqnarray}
\langle h_t^m,\varphi\rangle\to \langle h_t,\varphi\rangle,\ (\Delta_x h ^m, \, \Delta_x \varphi)\to (\Delta_x h,\, \Delta_x  \varphi)
  \label{convergence2}
\end{eqnarray}
as $m\to\infty$, for all $\varphi \in L^2(0,T;H_N^2(\Omega))$.

Taking the inner product of \eq{eq1Smoothed} and $ \varphi$  we arrive at
\begin{eqnarray}
  0 &=&  \langle\frac{\partial h^m}{\partial t}, \varphi\rangle + \nu(\Delta_x  h^m ,\, \Delta_x \varphi)
   -  (J(\widetilde{\nabla_x  h^m}), \nabla_x \varphi) \non\\
  &\to &  \langle\frac{\partial h}{\partial t}, \varphi\rangle + \nu(\Delta_x  h ,\, \Delta_x \varphi)
   -     (  J( {\nabla_x h}), \nabla_x \varphi).
  \label{eq1SmoothedWeak}
\end{eqnarray}
Thus \eq{definition} is proved. From \eq{est1a}, \eq{convergence0} and the choice of the smooth initial data $h_0^\varepsilon$ (let $\varepsilon=
\frac1m$) we have
\begin{eqnarray}
 \int_\Omega h^m(t,x) dx = \int_\Omega h_0^{\frac1m}(x) dx\to  \int_\Omega h_0 (x) dx,
\end{eqnarray}
and the left hand side converges to $\int_\Omega h (t,x) dx$, thus the mass is conserved for weak solution. And the existence of
weak solutions is complete.

\bigskip Next we are going to study the \\
\noindent{\bf Stability and Uniqueness.} Let $h_1,\ h_2$ be two weak solutions corresponding
to initial data $h_0^1$ and $h_0^2$, respectively. Define $u = h_1 - h_2$. We write
$$
 \nabla_x h_i = (p_i,q_i),\ (i=1,2),\  \nabla_x u = (p,q),\ J_i = J(\nabla_x h_i).
$$
Then by the estimates in Lemma~\ref{L2.0} we have $|J_1 - J_2 |\le C|\nabla_x h_1 - \nabla_x h_2| = C|\nabla_x u|  $, hence
\begin{eqnarray}
 0 &=& \frac12\frac{d}{dt}\|u\|^2 + \nu\|\Delta_x u\|^2 + \int_\Omega (J_1
 - J_2 )\cdot  \nabla_x u  dx \non\\
 &\ge & \frac12\frac{d}{dt}\|u\|^2 + \nu\|\Delta_x u\|^2 - C \|\nabla_x u\|^2.
 \label{unique1}
\end{eqnarray}

Making use of the Nirenberg inequality of the following form
\begin{eqnarray}
 \|\nabla_xu\| \le C\|\Delta_xu\|^\frac12 \| u\|^\frac12 + C' \| u\| ,
 \label{nirenbergL2}
\end{eqnarray}
and the Young inequality, from  \eq{unique1} one obtains
\begin{eqnarray}
  \frac12\frac{d}{dt}\|u\|^2 + \nu\|\Delta_x u\|^2  &\le &  C (\|\Delta_x u\| \|  u\| + \| u\|^2)   \non\\
 &\le & \frac\nu2\|\Delta_x u\|^2 + C \|u\|^2.
 \label{unique6}
\end{eqnarray}
Now using the Gronwall inequality we  get
\begin{eqnarray}
  \|u(t)\|^2   \le \| u(0)\|^2  e^{ C\, t} .
 \label{unique7}
\end{eqnarray}
Here $u(0) = h_0^1 - h_0^2$.   Thus the solution depends continuously on the initial data.

Consequently, if $h_0^1 = h_0^2$, that is $ \|u(t)\|^2 = 0$  which implies $\|u(t)\|^2 =0$, so the weak solution is unique.
 Therefore, the proof of Theorem~1.3 is complete.

\section{Regularity of weak solutions}

We shall investigate the regularity of weak solutions to both models with and without slope selection, while for the latter
model we can only carry out such study in one space dimension.

\subsection{The model with slope selection}

Suppose now that $h_0\in H^{2m}(\Omega)$ with $m\in\cN$ and there exists a unique solution $h$ to problem \eq{eq1} -- \eq{eq3}. In this section
we shall investigate the regularities of this solution.  We first consider the case that $m=1$.
\begin{Lemma}\label{L3.1} There hold for $h_0\in H^2(\Omega)$ that
\begin{eqnarray}
 \|h(t)\|^2_{H^2(\Omega)}  + \|h _t\|^2_{L^2(0,T;L^2(\Omega)}  &\le & C,
  \label{est3a}\\
   \|  h \|^2_{L^2(0,T;H^4(\Omega))}&\le & C.
  \label{est3b}
\end{eqnarray}

\end{Lemma}

\noindent{\it Proof.}  Multiplying equation \eq{eq1} by $h_t$ and integrating the resulting equation
with respect to $x$ yield
\begin{eqnarray}
 0 &=& \|h_t \|^2 + \frac\nu2\frac{d}{dt}  \|\Delta_x h \|^2 +\int_\Omega J\cdot \nabla_x h_t dx \non\\
  &=& \|h_t \|^2 + \frac{d}{dt} E[h] (t).
  \label{est3a1}
\end{eqnarray}
Thus one has
\begin{eqnarray}
\int_0^T\|h_t \|^2 dt&\le & C, \non\\
  E[h] (t) &\le & E[h] (t) .
  \non
\end{eqnarray}
Recalling the definitions of $E[h]$ and $f$ we arrive at \eq{est3a}.

From equation \eq{eq1}, making use the estimates in Lemma~\ref{L2.0} and \eq{est3a} we obtain
\begin{eqnarray}
 \int_0^T\|\Delta_x^2 h(\tau) \|^2 d\tau &\le & C\int_0^T\|h_t(\tau) \|^2 d\tau + C\int_0^T\|D^2 h\|^2 d\tau, \non\\
 &\le & C .  \non
 \end{eqnarray}
Thus \eq{est3b} is proved. And the proof of this lemma is complete.

\medskip
 To get the a priori estimates for higher order derivatives,
we differentiate the equation with respect to $t$ to get
\begin{eqnarray}
 \frac{\partial h_t}{\partial t}  + \nu\Delta_x ^2 h_t + {\rm div}_x \left( (J(\nabla_x h) )_t\right)  = 0 .
  \label{eq1t}
\end{eqnarray}
Such computations are formal. However by using the technique of difference quotient one can
justify easily. In a similar way for deriving the estimates for $h$, we arrive at
\begin{Lemma} \label{L3.2} Suppose that $h_0\in H^4(\Omega)$, i.e. $m=2$. There holds for any $t\in [0,T]$ that
\begin{eqnarray}
 \|h_t \|^2 + \|h \|^2_{H^4(\Omega)} + \|h_t\|^2_{L^2(0,T;H^2(\Omega)} \le C.
 \label{est1t}
\end{eqnarray}

\end{Lemma}

\noindent{\it Proof.}  Multiplying equation \eq{eq1t} by $h_t$ and integrating the resulting equation
with respect to $x$ yield
\begin{eqnarray}
 0 &=& \frac12\frac{d}{dt}  \|h_t \|^2 +  \nu \|\Delta_x h_t \|^2 +\int_\Omega J_t \cdot \nabla_x h_t dx \non\\
  &\ge & \frac12 \frac{d}{dt} \|h_t \|^2 +  \nu \|\Delta_x h_t \|^2 - C \|\nabla_x h_t \|^2 .
  \label{est3a1a}
\end{eqnarray}

With the help of the Nirenberg inequality \eq{nirenbergL2}
one obtains from  \eq{est3a1a} and the Young inequality that
\begin{eqnarray}
  \frac12 \frac{d}{dt} \|h_t \|^2 +  \nu \|\Delta_x h_t \|^2
  &\le & C  \|\Delta_xh_t \| \, \|h_t \| +    C \| h_t \|^2\non\\
  &\le & \frac\nu2  \|\Delta_xh_t \| ^2 +    C \| h_t \|^2 .
  \label{est3a2}
\end{eqnarray}
From which, by the Gronwall inequality,  it follows that
\begin{eqnarray}
 \|h_t \|^2 +  \nu \int_0^T\|\Delta_x h_t \|^2dt  &\le & C .
  \label{est3a3}
\end{eqnarray}
Furthermore, one can get from equation \eq{eq1} and Lemma~\ref{L2.0} that
\begin{eqnarray}
  \|\Delta^2_x h \|^2 &\le & C\|h_t \|^2 + C \|\Delta^2_x h \|^2\non\\
    &\le&  C,
  \label{est2e}
\end{eqnarray}
which implies that $h \in L^\infty(0,T;H^4(\Omega))$. The proof of this lemma is thus complete.

For the higher order derivatives with both $x$ and $t$,  we have
\begin{Lemma} \label{L3.3} Let $h_0\in H^{2m}(\Omega)$.
There holds for any $t\in [0,T]$ that
\begin{eqnarray}
 \|h (t)\|^2_{H^{2m}(\Omega)} + \|h \|^2_{L^2(0,T;H^{2m+2}(\Omega)} + \|h _t \|^2_{L^2(0,T;H^{2m-2}(\Omega)}  \le  C,
 \label{est2t1}
\end{eqnarray}
and
\begin{eqnarray}
 \|D^l_t h \|_{L^\infty(0,T;L^2(\Omega))}^2 + \|D^l_t h \|^2_{L^2(0,T;H^{2 }(\Omega)}  &\le& C, \ if\ m = 2l,
 \label{est2t2a}\\
  \|D^{l-1}_t h \|_{L^\infty(0,T;H^2(\Omega))}^2 + \|D^l_t h \|^2_{L^2(0,T;L^{2}(\Omega)} &\le&  C,\ if\ m = 2l-1 .
 \label{est2t2b}
\end{eqnarray}

\end{Lemma}

\noindent{\it Proof.} We employ the mathematical induction. From Lemma~3.1, it is easy to see that  \eq{est2t1} and
 \eq{est2t2b} are true for $m=1$ which implies $l=1$ too. By Lemma~3.2, estimate \eq{est2t2a} holds when $m=2$.
 Assume that  \eq{est2t1} is true for any $k\le m\in\cN$ and  \eq{est2t2a} and  \eq{est2t2b} are true respectively for even and odd $m$.
Next we shall prove they are true for $k\le m+1$ when $h_0\in H^{2(m+1)}(\Omega)$.

For the case that $m+1$ is even (resp. odd), differentiating $\frac{m+1}{2}$ (resp. $\frac{m}{2}$)
times equation \eq{eq1} with respect to $t$, letting $v=D^{\frac{m+1}{2}}_th$ (resp. $v=D^{\frac{m}{2}}_t h$), repeating the argument of Lemma~3.2 (resp. Lemma~3.1) for this function  $v$, and using the estimates in Lemma~2.1,
 we then conclude \eq{est2t1} holds for $m+1$, moreover, \eq{est2t2a} (resp. \eq{est2t2b}) is true for $l=\frac{m+1}{2}$ (resp. $l=\frac{m}{2}$).
  Thus  the proof of this lemma is complete.

\subsection{The one-dimensional cubic model without slope selection}

In this subsection we shall study the one-dimensional problem with a cubic current, i.e. \eq{nonlinearity1}, or \eq{nonlinearity2a}  -- \eq{nonlinearity2}, which is studied in  \cite{Stroscio95,Rost97,Watson06}. However the original two dimensional problem
 is still open. The problem is
\begin{eqnarray}
 h_t   + \nu h_{xxxx}  + (J(h_x))_x &=&0  \ {\rm in}\ Q_T,
  \label{cubic1}  \\
 h_x =0,\    h_{xxx}  &=&0  \ {\rm on}\ [0,T]\times \partial \Omega,
  \label{cubicBC}\\
  h|_{t=0}    &=& h_0.
  \label{cubicIC}
\end{eqnarray}
Here $\Omega=(a,b)\subset\cR$, $a,b\in\cR$ and  $J = \alpha h_x(1-h_x^2)$.

Assume that $h_0\in L^2(\Omega)$. Multiplying \eq{cubic1} by $h$ and integrating it with respect to $x$ give
\begin{eqnarray}
 0 &=& \frac12 \frac{d}{dt}\|h\|^2   + \nu \|h_{xx}\|^2  -\int_\Omega  J(h_x) h_x dx\non\\
  &=&\frac12 \frac{d}{dt}\|h\|^2   + \nu \|h_{xx}\|^2 + \int_\Omega    |h_x|^4 dx - \|h_{x}\|^2,
  \label{cubic2}
\end{eqnarray}
By the Young inequality $a^2\le \frac12a^4 + \frac12$, we infer from \eq{cubic2} that
\begin{eqnarray}
  \frac12 \frac{d}{dt}\|h\|^2   + \nu \|h_{xx}\|^2 + \int_\Omega    |h_x|^4 dx  \le  \frac12 \int_\Omega    |h_x|^4 dx + \frac12,
  \label{cubic3}
\end{eqnarray}
which gives

\begin{Lemma} \label{L3.4} We have
\begin{eqnarray}
 \|h  \|^2  + \int_0^T \left( \|h (t)\|^2_{H^{2}(\Omega)}  + \int_\Omega    |h_x|^4 dx\right)dt \le C.
 \label{estCubic1}
\end{eqnarray}

\end{Lemma}

Based on this lemma we can define weak solutions and prove the existence and uniqueness
 in a similar way as in Section~2 for $h_0\in L^2(\Omega)$. Suppose that
 $h_0\in H^1(\Omega)$, is the weak solution regular, i.e. $h \in L^2(0,T;H^1(\Omega))$? To answer this question, we need more
 estimates. From equation \eq{cubic1}, we obtain
\begin{eqnarray}
 0 &=& \frac12 \frac{d}{dt}\|h_x\|^2   + \nu \|h_{xxx}\|^2 - \int_\Omega J(h_x)_x h_{xx} dx,
  \label{cubic4}\\
  &=& \frac12 \frac{d}{dt}\|h_x\|^2   + \nu \|h_{xxx}\|^2 + \int_\Omega  3|h_x|^2  |h_{xx}|^2 dx - \|h_{xx}\|^2.
  \label{cubic4b}
\end{eqnarray}
Therefore,
$$
 \|h_x\|^2 + \int_0^T\left(\nu\|h_{xxx}\|^2 + \int_\Omega  3|h_x|^2 |h_{xx}|^2 dx\right)dt\le C.
$$
So the weak solution is more regular if $h_0\in H^1(\Omega)$.

However for two-dimensional problem, we can not get \eq{cubic4} so that we obtain the good term
 $\int_\Omega  3|h_x|^2  |h_{xx}|^2 dx $ in \eq{cubic4b}.

\bigskip
\bigskip
\noindent{\bf Acknowledgement.} The first author (JZ) of this work has
been partly supported National Natural Science Foundation of China (Grant No. 10501040)
and by Zhejiang Provinicial Natural Science Foundation of China (Grant No. Y6100611),
 and the second (PZ) by Grant MTM2008-03541 of
the Ministerio de Educac\'ion y Ciencia of Spain, and by Project
PI2010-04 of the Basque Government. The second author would like to
express his sincere thanks to Dr. S. Watson at Univ. of Glasgow for
 bringing the paper \cite{Watson06} to the authors' attention.

\end{document}